\documentclass{amsart}

\usepackage[cmtip,all]{xy}

\usepackage{amsmath,amscd,amssymb,hyperref,stmaryrd,verbatim,mathtools}

\copyrightinfo{2017}{American Mathematical Society}

\newtheorem{theorem}{Theorem}[section]

\newtheorem{prop}[theorem]{Proposition}
\newtheorem{cor}[theorem]{Corollary}

\theoremstyle{definition}

\newtheorem{example}[theorem]{Example}

\theoremstyle{remark}
\newtheorem{remark}[theorem]{Remark}

\numberwithin{equation}{section}

\DeclareMathOperator{\tc}{{\sf TC}}
\DeclareMathOperator{\zcl}{\sf zcl}
\DeclareMathOperator{\secat}{secat}
\DeclareMathOperator{\cat}{cat}
\DeclareMathOperator{\cd}{cd}
\newcommand{\A}{{\mathcal{A}}}
\newcommand{\C}{{\mathbb{C}}}
\newcommand{\R}{{\mathbb{R}}}
\newcommand{\Z}{{\mathbb{Z}}}
\newcommand{\bx}{{\mathbf{x}}}
\DeclareMathOperator{\PSL}{PSL}
\DeclareMathOperator{\SO}{SO}
\DeclareMathOperator{\hdim}{hdim}

\begin{document}

\title[Topological complexity of configuration spaces]{Topological complexity of classical configuration spaces and related objects}

\author{Daniel C. Cohen}
\address{Department of Mathematics, Louisiana State University, Baton Rouge, LA, USA}
\email{\href{mailto:cohen@math.lsu.edu}{cohen@math.lsu.edu}}
\urladdr{\href{http://www.math.lsu.edu/~cohen/}{www.math.lsu.edu/\char'176cohen}}
\thanks{Partially supported by NSF 1105439}



\begin{abstract}
We survey results on the topological complexity of classical configuration spaces of distinct ordered points in orientable surfaces and related spaces, including certain orbit configuration spaces and Eilenberg-Mac\,Lane spaces associated to certain discrete groups.
\end{abstract}

\maketitle

\setcounter{tocdepth}{2}
\tableofcontents

\section{Introduction} \label{sec:intro}
Investigation of the collision-free motion of $n$ distinct ordered particles in a topological space $X$ leads one to study the (classical) configuration space
\[
F(X,n)=\{(x_1,\dots,x_n) \in X^n \mid x_i \neq x_j\ \text{if}\ i\neq j\}
\]
of $n$ distinct ordered points in $X$, 
and the topological complexity of this space. For a path-connected topological space $Y$, $I=[0,1]$ the unit interval, and $Y^I$ the space of all continuous paths $\gamma\colon I \to Y$ (with the compact-open topology), the \emph{topological complexity} of $Y$ is the sectional category (or Schwarz genus) of the fibration $\eta\colon Y^I \to Y\times Y$, $\gamma \mapsto (\gamma(0),\gamma(1))$, $\tc(Y)=\secat(\pi)$. This homotopy invariant, introduced by Farber \cite{FaTC1}, provides a topological approach to the motion planning problem from robotics.

For any $s\ge 2$, one may more generally consider the sectional category of the fibration $\eta_s\colon Y^I \to Y^{s}$, where $\eta_s$ sends a path to $s$ points on the path, $\eta_s(\gamma)=(\gamma(t_1),\dots,\gamma(t_s))$, where $t_k=(k-1)/(s-1)$, $1\le k\le s$. This is the higher topological complexity of $Y$, $\tc_s(Y)=\secat(\eta_s)$ introduced by Rudyak \cite{Ru}, extending Farber's notion above, as $\tc_2(Y)=\tc(Y)$.

We survey results on the topological complexity of configuration spaces $F(X,n)$ in the case where $X$ is an orientable surface, as well as related objects. 
The discussion is focused primarily on the ``classical'' topological complexity, 
$\tc=\tc_2$, 
and includes remarks on the higher topological complexity in those instances where results on this invariant are known.  
The general principle is as follows:

\smallskip

\emph{The topological complexity is as large as possible, given natural constraints.}

\smallskip

\noindent For instance, as is well known (and discussed in Section \ref{sec:plane} below), the configuration space of $n$ distinct ordered points in the plane has the homotopy type of the product of a circle and a CW-complex of dimension $n-2$. These constraints, together with the fact that $\tc(S^1)=2$, and the known bounds recorded next, yield the topological complexity of $F(\R^2,n)$, recorded in Theorem \ref{thm:plane}.

If $Y$ is a topological space with the homotopy type of a finite-dimensional CW-complex, let $\mathrm{hdim}(Y)$ denote the homotopy dimension of $Y$. Throughout the discussion, we will make use of the following basic tools. For details and other relevant facts, see Farber's survey \cite{dc_Fa05}. 
\[
\begin{aligned}
&\bullet \quad\tc(Y) \le 2\cdot \mathrm{hdim}(Y)+1 \qquad  \bullet\quad \tc(Y\times Z) \le \tc(Y)+\tc(Z)-1\\
&\bullet \quad\tc(Y) > \zcl H^*(Y) = \mathrm{cup\ length} \bigl[\ker \bigl(H^*(Y) \otimes H^*(Y) \xrightarrow{\ \cup\ } H^*(Y)\bigr)\bigr]
\end{aligned}
\]
We call the first two of these the dimension and product inequalities, and use cohomology with $\C$-coefficients (unless stated otherwise) in the context of the third, the zero divisor cup length. 
We use the unreduced notions of topological complexity and higher topological complexity. For instance, $\tc(Y)=\secat(\eta\colon Y^I \to Y\times Y)$ is equal to the smallest integer $m$ such that there exists of cover of $Y\times Y$ by $m$ open sets, on each of which the fibration $\eta\colon Y^I \to Y\times Y$ admits a continuous local section. In particular, for $Y$ contractible, $\tc(Y)=1$.

\section{The plane and the sphere} \label{sec:plane}
The topological complexity of the configuration space of ordered points in the plane $X=\R^2=\C$ was determined by Farber-Yuzvinsky.
\begin{theorem}[\cite{dc_FY}] 
\label{thm:plane}
$\tc(F(\C,n))=2n-2$ for $n\ge 2$.
\end{theorem}

\subsection{Arrangements, I} \label{subsec:arr1}

To discuss this result, we recall some relevant facts from the theory of complex hyperplane arrangements. 
Let $V=\C^\ell$ be a complex vector space of dimension $\ell<\infty$. 
A hyperplane in $V$ is a codimension one affine subspace $H$. A hyperplane arrangement in $V$ is a finite collection $\A=\{H_1,\dots,H_m\}$ of hyperplanes in $V$. 
If we fix coordinates $\bx=(x_1,\dots,x_\ell)$ on $V$, each hyperplane of $\A$ may be realized as $H_j=\{\bx \in \C^\ell\mid f_j(\bx)=0\}$, where $f_j$ is a linear polynomial.  The product $Q(\A)=\prod_{j=1}^m f_j$ is said to be a defining polynomial of the arrangement~$\A$. 

An arrangement $\A$ is said to be essential if there are $\ell$ hyperplanes in $\A$ whose intersection is a point. If the intersection of all hyperplanes in $\A$ is nonempty, $\bigcap_{j=1}^m H_j \neq \emptyset$, we can choose coordinates so that $f_j$ is a linear form and $0 \in H_j$ for each $j$. In this situation, $\A$ is said to be central, and the defining polynomial $Q(\A)$ is homogeneous. Refer to Orlik-Terao \cite{OT} as a general reference on arrangements.

A principal object of topological study of arrangements is the complement $M=M(\A)=V\smallsetminus \bigcup_{j=1}^m H_j = V\smallsetminus Q(\A)^{-1}(0)$. The complement is an open, smooth manifold of real dimension $2\ell$, which has the homotopy type of a connected, finite CW-complex of dimension at most $\ell$ ($M$ is a Stein manifold). If $\A$ is essential, then the homotopy dimension of $M$ is precisely $\ell$.

\begin{example} \label{ex:braid arr}
The braid arrangement $\A_n$ in $V=\C^n$ is the central arrangement consisting of the $\binom{n}{2}$ diagonal hyperplanes $H_{i,j}=\{\bx\in\C^n \mid x_i-x_j=0\}$. Note that $\A_n$ is not essential, as the diagonal line $\{(x,x,\dots,x)\}$ is contained in (the intersection of) all hyperplanes of $\A_n$. 
The complement of the braid arrangement
\[
M(\A_n)=\C^n \smallsetminus \bigcup_{i<j} H_{i,j}=\{(x_1,\dots,x_n) \in \C^n \mid x_i \neq x_j\ \text{if}\ i \neq j\}=F(\C,n)\\
\]
is the configuration space of $n$ distinct ordered points in the plane $\C$, with fundamental group $\pi_1(F(\C,n))=P_n$, the Artin pure braid group on $n$ strings.

Projecting along the diagonal line $x_1=\dots=x_n$ onto the subspace $\{\bx \in \C^n \mid \sum x_i=0\}\cong \C^{n-1}$ yields an essential arrangement $\hat \A_n$ in $\C^{n-1}$. Denoting the coordinates of $\C^{n-1}$ by $\mathbf{y}=(y_1,\dots,y_{n-1})$ (where $y_i=x_i-x_n$), the hyperplanes of $\hat \A_n$ are $\{y_i=0\}$, $1\le i\le n$, and $\{y_i-y_j=0\}$, $1\le i<j\le n-1$. Letting $M_n=M(\hat\A_n)$ be the complement, we have $F(\C,n)=M(\A_n) \cong M_n \times \C$, and $F(\C,n)\simeq M_n$ has the homotopy type of a CW-complex of dimension $n-1$.  Note that $M_n = F(\C^*,n-1)$ is itself a configuration space.
\end{example}

The diagonal action of the group $\C^*$ on $\C^{\ell}\smallsetminus\{0\}$, $\xi \cdot (x_1,\dots,x_\ell)=
(\xi x_1,\dots,\xi x_\ell)$, is free. The orbit space and map are the complex projective space ${\C}P^{\ell-1}$ and the Hopf fibration $p\colon \C^\ell\smallsetminus\{0\} \to {\C}P^{\ell-1}$, respectively.
For any nonempty central arrangement $\A$ in $\C^\ell$, with complement $M=M(\A)$, the restriction of the Hopf map, $p\colon M \to p(M)$, is a bundle map, with fiber $\C^*$. If $H_1\in\A$, we have $p(H_1)\cong {\C}P^{\ell-2}$ and ${\C}P^{\ell-1}\smallsetminus p(H_1) \cong \C^{\ell-1}$ is contractible. Thus, $p\colon M \to p(M)$ is the restriction of a trivial bundle.
\[
\begin{CD} \C^* @>>> \C^* @>>> \C^*\\
@VVV  @VVV  @VVV\\
M @>>> \C^\ell \smallsetminus H_1 @>>>  \C^\ell\smallsetminus\{0\} \\
@VVV @VVV @VVpV\\
p(M) @>>> {\C}P^{\ell-1}\smallsetminus p(H_1) @>>>{\C}P^{\ell-1}
\end{CD}
\]
\begin{prop} \label{prop:central/affine}
Let $\A$ be a nonempty central arrangement in $\C^\ell$. Then the complement $M$ of $\A$ is diffeomorphic to the product space $p(M) \times \C^*$.
\end{prop}

\begin{remark} \label{rem:decone}
The projective image $p(M)\subset{\C}P^{\ell-1}$ of the complement of the central arrangement $\A$ of $n\ge 1$ hyperplanes in $\C^\ell$ may itself be realized as the the complement of a (not necessarily central) arrangement of hyperplanes in the affine space $\C^{\ell-1}$. If the hyperplanes $H_j$ of $\A$ are defined by linear forms $f_j$, we may choose coordinates $x_1,x_2,\dots,x_\ell$ on $\C^\ell$ so that $f_1=x_1$, that is, $H_1=\{x_1=0\}$. Then, $p(M)\cong M(\mathrm{d}\A)$, where $\mathrm{d}\A$ is decone of $\A$ with respect to $H_1$, the arrangement in $\C^{\ell-1}$, with coordinates $x_2,\dots,x_\ell$, defined by the linear polynomials $f_j(1,x_2,\dots,x_\ell)$, $2 \le j \le n$.
\end{remark}

\subsection{Cohomology} \label{subsec:cohomology}
The cohomology of the configuration space $F(X,n)$ of $n$ distinct ordered points in $X$ has been the object of a great deal of study, particularly for $X$ a manifold.  See, for instance, \cite{dc_To} and the references therein. In the case where $X$ is a Euclidean space, the structure of the ring $H^*(F(\R^m,n);\Z)$ was determined by Arnol'd \cite{dc_Ar} (for $\R^2=\C$) and Cohen \cite{dc_Co}.  
\begin{theorem} 
\label{thm:configcoho}
The integral cohomology ring $H^*(F(\R^m,n);\Z)$ has generators $\alpha_{i,j}^{}$, $1\le i,j\le n$, $i\neq j$, of degree $m-1$, and relations $\alpha_{i,j}^{}=(-1)^m\alpha_{j,i}^{}$, $\alpha_{i,j}^2=0$, and
$\alpha_{i,j}^{}\alpha_{i,k}^{}+\alpha_{j,k}^{}\alpha_{j,i}^{}+\alpha_{k,i}^{}\alpha_{k,j}^{}=0$
for distinct $i,j,k$. 
\end{theorem}

The cohomology of the complement of an arbitrary complex hyperplane arrangement was determined by Brieskorn \cite{dc_Br} (resolving positively a conjecture of Arnol'd).
\begin{theorem}
Let $\A=\{H_1,\dots,H_m\}$ be an arrangement in $\C^\ell$, where $H_j=\{f_j=0\}$. Then the integral cohomology of the complement $M$ is torsion free, and is generated by the cohomology classes of the $1$-forms 
$\frac{1}{2\pi i}d\log f_j$.
\end{theorem}
\begin{remark}
Let $R$ denote the algebra of differential forms generated by $1$ and the $1$-forms $\alpha_j=\frac{1}{2\pi i}d\log f_j$, $1\le j \le m$. Brieskorn's theorem shows that the inclusion of $R$ in the 
de\,Rham complex of smooth forms on $M$ is a quasi-isomorphism, as $\alpha_j \mapsto [\alpha_j]$ induces an isomorphism $R \to H^*(M;\R)$. Consequently, the complement of a complex hyperplane arrangement is a (rationally) formal space.
\end{remark}

We now discuss how the above considerations may be used to establish Theorem \ref{thm:plane}, determining the topological complexity of the configuration space $F(\C,n)$.  From Example \ref{ex:braid arr}, we have $F(\C,n) \simeq M_n$, and by Proposition \ref{prop:central/affine}, $M_n \cong p(M_n) \times \C^*$. Since $p(M_n)$ has the homotopy type of a CW-complex of dimension $n-2$, the product inequality yields 
\[
\tc(F(\C,n))=\tc(M_n) \le \tc(p(M_n))+\tc(\C^*)-1 \le 2n-2.
\]
The reverse inequality $\tc(F(\C,n)) \ge 2n-2$ is obtained using the zero divisor cup length. Let $A=H^*(F(\C,n);\C)$ and, abusing notation, denote the generators of $A$ (from Theorem \ref{thm:configcoho} and the Universal Coefficient Theorem) by $\alpha_{i,j}$, $1\le i<j\le n$.

\begin{prop}[\cite{dc_FY}] The zero-divisors $\bar{\alpha}_{i,j}=1\otimes \alpha_{i,j}-\alpha_{i,j}\otimes 1 \in A^1 \otimes A^1$ satisfy $\bar{\alpha}_{1,2}\cdot\bar{\alpha}_{1,3} \cdots \bar{\alpha}_{1,n} \cdot \bar{\alpha}_{2,3} \cdots \bar{\alpha}_{2,n} \neq 0$. Consequently, $\zcl H^*(F(\C,n)) \ge 2n-3$.
\end{prop}
With the above considerations, this yields $\tc(F(\C,n))=2n-2$. Similarly, the topological complexity of the configuration space of ordered points in the punctured plane was determined by Farber-Grant-Yuzvinsky.
\begin{theorem}[\cite{dc_fgy}] \label{thm:fgy}$\tc(F(\C\smallsetminus \{m\ \text{points}\},n))=
\begin{cases} 2n,&\text{if $m=1$.}\\ 2n+1,&\text{if $m\ge 2$.}
\end{cases}$
\end{theorem}

For $m=1$, since $F(\C^*,n) \simeq F(\C,n+1)$, this is a restatement of Theorem \ref{thm:plane}. For $m\ge 2$, since 
$F(\C\smallsetminus\{m\ \text{points}\},n)$ may be realized as the complement of an essential arrangement in $\C^n$, we have $\hdim F(\C\smallsetminus\{m\ \text{points}\},n)=n$, and consequently $\tc(F(\C\smallsetminus\{m\ \text{points}\},n)) \le 2n+1$. The reverse inequality is obtained by showing that the zero divisor cup length of the cohomology ring is (at least) $2n$. See Section \ref{sec:groups} for additional complementary discussion.

\begin{remark} The topological complexity of the configuration space of ordered points in a higher dimensional Euclidean space is also known:
\[
\tc(F(\R^k,n))=\begin{cases} 2n-1,&\text{if}\ k\ge 3 \ \text{odd\hskip 6pt \cite{dc_FY},}\\ 2n-2,&\text{if}\ k\ge 4\ \text{even \cite{dc_FG}}.
\end{cases}
\]
\end{remark}

The higher topological complexity of the configuration space of ordered points in a higher dimensional (punctured) Euclidean space was computed by Gonz\'alez-Grant.
\begin{theorem}[\cite{dc_GoGr}] \label{thm:gogr}
Let $k,m,n,s$ be nonnegative integers with $k,s\ge 2$, $n\ge 1$, and such that $n\ge 2$ if $m=0$. Then
\[
\tc_s(F(\R^k\smallsetminus \{m\ \text{points}\},n))=
\begin{cases} s(n-1),&\text{if $m=0$ and $k \equiv 0\ \mathrm{mod}\ 2$;}\\ 
s(n-1)+1,&\text{if $m= 0$ and $k \equiv 1\ \mathrm{mod}\ 2$;}\\
sn, &\text{if $m=1$ and $k \equiv 0\ \mathrm{mod}\ 2$;}\\
sn+1,&\text{otherwise.}
\end{cases}
\]
\end{theorem}

\subsection{Genus zero}\label{sec:sphere}
The topological complexity of the configuration space of ordered points in the sphere $X=S^2$ was determined by Cohen-Farber, using results of Farber, Grant, and Yuzvinsky \cite{dc_Fa05,dc_fgy}.
\begin{theorem}[\cite{dc_CF}] $\tc(F(S^2,n))=
\begin{cases}3,&\text{if}\ n=1,2,\\2n-2,&\text{if}\ n\ge 3.
\end{cases}$
\end{theorem}

For $n\le 2$, the configuration space $F(S^2,n)$ has the homotopy type of $S^2$ itself, and it is well known that $\tc(S^2)=3$.

For $(x_1,x_2,x_3) \in F(S^2,3)$, there is a unique M\"obius transformation $A$ of $S^2={\C}P^1$ taking $x_1$ to $0$, $x_2$ to $1$  and $x_3$ to $\infty$, and depending continuously on $(x_1,x_2,x_3)$. This yields a homeomorphism from $F(S^2,3)$ to $\mathrm{PSL}(2,\C)$, the group of M\"obius transformations acting on the sphere, and for $n\ge 4$, a homeomorphism
\[
F(S^2, n)\to \PSL(2, \C) \times F(S^2\smallsetminus\{0,1,\infty\}, n-3),
\]
taking $(x_1, \dots, x_n)$ to $(A, (y_4, \dots, y_{n}))$. 
Here, we  have identified $S^2={\C}P^1$, 
$A$ is the above transformation, and $y_i=Ax_i$ for $i=4, \dots, n$.
Since $\PSL(2,\C)$ deformation retracts onto $\SO(3)$ and  $F(S^2\smallsetminus \{0,1,\infty\},n-3) =F(\R^2\smallsetminus\{2\ \text{points}\}, n-3)$, we obtain homotopy equivalences
$F(S^2,3)\simeq \SO(3)$, and, for $n\ge 4$,
\begin{equation} \label{eq:sphere htpy}
F(S^2,n) \simeq \SO(3) \times F(\R^2\smallsetminus\{2\ \text{points}\},n-3).
\end{equation}

Consequently, for $n=3$, we have $\tc(F(S^2,3))=\tc(\SO(3))=\cat(\SO(3))=4$, since $\SO(3)$ is a connected Lie group, see \cite[Lemma 8.2]{dc_Fa03}.

For $n\ge 4$, the product inequality yields
\[
\tc(F(S^2,n)) \le \tc(\SO(3))+\tc(F(\R^2\smallsetminus \{2\ \text{points}\},n-3))-1 = 2n-2,
\]
using the Farber-Grant-Yuzvinsky result concerning the topological complexity of $F(\R^2\smallsetminus \{2\ \text{points}\},n-3)=F(\C\smallsetminus \{2\ \text{points}\},n-3)$ recorded in Theorem \ref{thm:fgy} above. The reverse inequality $\tc(F(S^2,n)) \ge 2n-2$ may be obtained by showing that the zero divisor cup length of the cohomology ring 
\[
H^*(F(S^2,n);\Z_2) \cong H^*(\SO(3);\Z_2) \otimes H^*(F(\C\smallsetminus \{2\ \text{points}\},n-3);\Z_2)
\] 
is (at least) $2n-3$.

Similar considerations yield the higher topological complexity of the configuration space of ordered points in the sphere, as determined by Gonz\'alez-Guti\'errez.

\begin{theorem}[\cite{dc_GG}]
For $s\ge 2$, $\tc_s(F(S^2,n))=\begin{cases}
s+1,&\text{if $n=1, 2$,}\\
sn-2,&\text{if $n\ge 3$.}
\end{cases}$
\end{theorem}

\section{Genus one}\label{sec:torus}
The topological complexity of the configuration space of ordered points in the torus $X=T=S^1\times S^1$ was determined by Cohen-Farber.
\begin{theorem}[\cite{dc_CF}] 
\label{thm:torus}
$\tc(F(T,n))=2n+1$.
\end{theorem}

For $n=1$, we have $F(T,1)=T$, and $\tc(T)=\tc(S^1\times S^1)=3$.

For $n\ge 2$, using the fact that $T$ is a group, it is readily checked that the map
\[
\bigl((u,v),\bigl((z_1,w_1),\dots,(z_{n-1},w_{n-1})\bigr)\bigr) \mapsto \bigl((u,v),(uz_1,vw_1),\dots,(uz_{n-1},vw_{n-1})\bigr)
\]
is a homeomorphism $T \times F(T\smallsetminus\{1 \ \text{point}\},n-1) \to F(T,n)$.

\subsection{Fadell-Neuwirth theorem}
The homotopy dimension of the configuration space $F(T\smallsetminus\{1 \ \text{point}\},n-1)$ may be obtained using the classical Fadell-Neuwirth theorem.

\begin{theorem}[\cite{dc_FN}] 
\label{thm:FN}
Let $X$ be a manifold without boundary of dimension at least two. For $\ell \le n$, the projection onto the first $\ell$ coordinates, 
$p:F(X,n)\to F(X,\ell)$, is a locally trivial bundle, with fiber 
$F(X\setminus\{\ell\ \text{points}\},n-\ell)$.
\end{theorem}
\begin{remark}
\label{rem:FNsection}
The aforementioned Fadell-Neuwirth bundles often admit cross sections. This is the case, for instance, 
\begin{enumerate}
\item[(i)] for the bundle $F(X,n) \to F(X,\ell)$, when $X$ is a punctured surface; and
\item[(ii)] for the bundle $F(X,n) \to F(X,1)=X$, when $X$ is a compact manifold with nonvanishing first Betti number.
\end{enumerate}
\end{remark}

Let $G=\pi_1(F(T\smallsetminus\{1 \ \text{point}\},n-1))$. Taking $X=T\smallsetminus\{1 \ \text{point}\}$ in the Fadell-Neuwirth theorem, the bundle $F(X,k) \to F(X,k-1)$ admits a cross section, and has fiber $X \smallsetminus\{k-1\ \text{points}\}=T\smallsetminus\{k\ \text{points}\}$ homotopy equivalent to a bouquet of circles. Using this repeatedly for 
$k=n,n-1,n-2\ldots$ reveals that $F(T\smallsetminus\{1 \ \text{point}\},n-1)$ is a $K(G,1)$-space, and that $G$ is an iterated semidirect product of (finitely generated) free groups. This being the case, the cohomological and geometric dimensions of $G$ are both equal to $n-1$, 
$\cd(G)=\mathrm{gd}(G)=n-1$, see \cite{dc_CS} and Section \ref{sec:groups} below. Consequently, $\mathrm{hdim}(F(T\smallsetminus\{1 \ \text{point}\},n-1))=n-1$. Then the dimension and product inequalities yield  $\tc(F(T,n))\le 2n+1$.

\subsection{Cohen-Taylor/Totaro spectral sequence}
The proof of Theorem \ref{thm:torus} is completed by showing that $\zcl H^*(F(T,n)) \ge 2n$. The tool here is the Cohen-Taylor/Totaro spectral sequence \cite{dc_CT,dc_To}. For $X$ a closed $m$-manifold, let $p_i\colon X^n \to X$ and $p_{i,j}\colon X^n \to X^2$ be the obvious projections. The inclusion $F(X,n) \to X^n$ yields a Leray spectral sequence converging to $H^*(F(X,n))$. The initial term is the quotient of the algebra $H^*(X^n) \otimes H^*(F(\R^m,n))$ by the relations $(p_i^*(u)-p_j^*(u)) \otimes \alpha_{i,j}$ for $i\neq j$, $u\in H^*(X)$, and $\alpha_{i,j}$ the generators of $H^*(F(\R^m,n))$ from Theorem \ref{thm:configcoho}. The first nontrivial differential is given by $d(\alpha_{i,j})=p_{i,j}^*(\Delta)$, where $\Delta \in H^m(X\times X)$ is the cohomology class dual to the diagonal. 
Explicitly, if $\Omega\in H^m(X)$ is a generator, and $\{\beta_i\}$ and $\{\beta_i^*\}$ are bases for $H^*(X)$ with $\beta_i \cup \beta_j^*=\delta_{i,,j}\Omega$, then $\Delta=\sum (-1)^{|\beta_i|} \beta_i \times \beta_i^*$, where $|\beta_i|$ is the degree of $\beta_i$, and $\delta_{i,j}$ is the Kronecker symbol.

As shown by Totaro \cite{dc_To}, for $X$ a smooth complex projective variety, this spectral sequence degenerates immediately, $d$ above is the {only} nontrivial differential.
\begin{prop}[\cite{dc_To,dc_CF}] \label{prop:subalg}
For $X$ a smooth complex projective variety, let $H=H^*(X^{n})$, and let $I$  be the ideal in $H$ generated by $\{p^*_{i,j}(\Delta) \mid i<j\}$. Then $H/I$ is a subalgebra of $H^*(F(X,n))$. Consequently, $\zcl H^*(F(X,n)) \ge \zcl H/I$ and $\tc(F(X,n)) \ge \zcl H/I+1$.
\end{prop}

For $X=T$, these considerations may be used to obtain the required lower bound on $\zcl H^*(F(T,n))$. Since $T^n=(S^1)^{2n}$, the algebra $H$ is graded commutative with degree-one generators. Denoting these generators by $a_i,b_i$, $1\le i \le n$, the ideal $I$ in $H$ is generated by the elements $p^*_{i,j}(\Delta)=(a_j-a_i)(b_j-b_i)$. 
The change of variables $x_1=a_1$, $y_1=b_1$, $x_j=a_j-a_1$, $y_j=b_j-b_1$ for $2\le j \le n$ 
reveals that the quotient $A=H/I$ 
is generated by degree-one classes $x_i$, $y_i$, $1\le i\le n$, with relations $x_iy_i=0$, $2\le i\le n$, $x_jy_k+x_ky_j=0$, $2\le j<k\le n$, and their consequences. Recall that, for an element $u\in A$, we denote by $\bar{u}=1\otimes u-u\otimes 1$ the corresponding zero divisor in $A\otimes A$.

\begin{prop}[\cite{dc_CF}] The zero divisors $\bar{x}_i$ and $\bar{y}_i$, $1\le i \le n$, in $A^1 \otimes A^1$ satisfy
$\bar{x}_1\cdot\bar{y}_1\cdot\bar{x}_2\cdot\bar{y}_2 \cdots \bar{x}_n\cdot\bar{y}_n \neq 0$. Consequently, $\zcl H^*(F(T,n)) \ge 2n$.
\end{prop}
With the above considerations, this yields $\tc(F(T,n))=2n+1$.

The higher topological complexity of the configuration space of the torus was determined by Gonz\'alez-Guti\'errez.

\begin{theorem}[\cite{dc_GG}] For $s\ge 2$, $\tc_s(F(T,n))=s(n+1)-1$.
\end{theorem}

\section{Higher genus}
The topological complexity of the configuration space of ordered points in a higher genus surface $X=\Sigma_g$, $g\ge 2$, was studied by Cohen-Farber, and Gonz\'alez-Guti\'errez.
\begin{theorem}[\cite{dc_CF,dc_GG}] \label{thm:high g}
$\tc(F(\Sigma_g,n))=2n+3$
\end{theorem}
The results of Gonz\'alez-Guti\'errez yield the higher topological complexity of the configuration space of a high genus surface.
\begin{theorem}[\cite{dc_GG}]
For $s\ge 2$, $\tc_s(F(\Sigma_g,n))=s(n+1)+1$.
\end{theorem}

The discussion below focuses on the classical topological complexity $\tc=\tc_2$. 
For $n=1$, $F(\Sigma_g,1)=\Sigma_g$, and $\tc(\Sigma_g)=5$, as is well known, see \cite[Theorem 9]{FaTC1}.

For $n\ge 2$, the configuration space $F(\Sigma_g,n)$ is an Eilenberg-Mac\,Lane space of type
$K(G,1)$, where $G=\pi_1(F(\Sigma_g,n))$ is the pure braid group of the surface $\Sigma_g$.  As noted in Remark \ref{rem:FNsection}\,(ii), the Fadell--Neuwirth fibration $F(\Sigma_g,n)\to\Sigma_g$ admits a section. Consequently, the surface pure braid group 
\[
G\cong \pi_1(F(\Sigma_g\smallsetminus \{1\ \text{point}\},n-1)) \rtimes \pi_1(\Sigma)
\] 
is a semidirect product. As in the genus one case, repeated application of the Fadell-Neuwirth theorem shows that the group $\pi_1(F(\Sigma_g\smallsetminus  \{1\ \text{point}\},n-1))$ is an $(n-1)$-fold iterated semidirect product of free groups.  It follows that the cohomological dimension of $G$ is equal to $n+1$, as is the geometric dimension.  Consequently, $F(\Sigma_g,n)$ has the homotopy type of a cell complex of dimension $n+1$.  So $\tc(F(\Sigma_g,n)) \le 2n+3$.

To complete the proof of Theorem \ref{thm:high g}, it remains to show that the zero divisors cup length is sufficiently large, $\zcl H^*(F(\Sigma_g,n)) \ge 2n+2$. 
As discussed in \cite{dc_GG}, the proof of this last fact given in \cite{dc_CF} contains an oversimplification, invalidating the argument there. A detailed argument establishing $\zcl H^*(F(\Sigma_g,n)) \ge 2n+2$ may be found in \cite{dc_GG}. 
We put forward below an alternate approach to this argument.

The cohomology ring $H^*(\Sigma_g)=H^*(\Sigma_g;\C)$ is generated by degree-one elements $a_k,b_k\in H^1(\Sigma_g)$, $1\le k \le g$, with relations $a_ra_s=a_rb_s=b_rb_s=0$ for $1\le r\neq s\le g$, and $a_kb_k=a_gb_g=\omega$ (a generator of $H^2(\Sigma_g)$) for $1\le k\le g-1$. The cohomology class $\Delta\in H^2(\Sigma_g\times\Sigma_g)$ dual to the diagonal may be expressed as
\[
\Delta = \omega\times 1+1\times \omega +\sum_{p=1}^g b_p\times a_p-a_p\times b_p.
\]
Let $H=H^*(\Sigma_g^n)=[H^*(\Sigma_g)]^{\otimes n}$, with generators $a_{i,p},b_{i,p}$, $1\le i\le n$, $1\le p \le g$, where, for instance,  $a_{i,p}=1\times \dots \times
\overset{i}{a_p} \times \dots \times 1$, and relations $a_{i,r}a_{i,s}=a_{i,r}b_{i,s}=b_{i,r}b_{i,s}=0$ for $1\le r\neq s\le g$ and $a_{i,k}b_{i,k}=a_{i,g}b_{i,g}=\omega_i$ for $1\le k\le g-1$, all for each $i$, $1\le i\le n$. In this notation, for the projection $p_{i,j}\colon \Sigma_g^n \to \Sigma_g^2$, we have
\[
p_{i,j}^*(\Delta)=\omega_i+\omega_j+\sum_{p=1}^g b_{i,p} a_{j,p}-a_{i,p} b_{j,p}
=a_{i,g}b_{i,g}+a_{j,g}b_{j,g}+\sum_{p=1}^g b_{i,p} a_{j,p}-a_{i,p} b_{j,p}.
\]
If $I=\langle p_{i,j}^*(\Delta) \mid 1\le i<j\le n\rangle$ is the ideal in $H$ generated by these elements, by Proposition \ref{prop:subalg} it suffices to show that $\zcl H/I \ge 2n+2$.

\subsection{Gr\"obner bases} \label{subsec:gb}
Write $A=H/I$. This algebra may be realized as the quotient of an exterior algebra by a homogeneous two-sided ideal. Let $E$ be the exterior algebra (over $\C$) generated by one-dimensional classes $a_{i,p},b_{i,p}$, $1\le i\le n$, $1\le p \le g$, and let $J$ be the ideal in $E$ given by
\[
J=\left\langle 
\begin{matrix}\ a_{i,r}a_{i,s},\ a_{i,r}b_{i,s},\ b_{i,r}b_{i,s},&(1\le i\le n, 1\le r\neq s\le g),\\ a_{i,k}b_{i,k}-a_{i,g}b_{i,g},&(1\le i\le n, 1\le k\le g-1),\\
\Delta_{i,j},&(1\le i<j\le n)\hfill
\end{matrix}\right\rangle,
\]
where $\Delta_{i,j}=a_{i,g}b_{i,g}+a_{j,g}b_{j,g}+\sum_{p=1}^g b_{i,p} a_{j,p}-a_{i,p} b_{j,p}$ in $E$.
Then $A=E/J$, and the ideal $J$ and the algebra $A$ may be studied using Gr\"obner basis theory in the exterior algebra, following Aramova-Herzog-Hibi \cite{AHH}. As in \cite{dc_GG}, it is convenient to work with a quotient of $A$. Let $K$ be the ideal in $E$ given by
\[
K=\langle a_{i,p}a_{j,q},\ a_{i,p}b_{j,q},\ b_{i,p}b_{j,q},\ (1\le i \neq j\le n, 2\le p,q \le g)\rangle,
\]
write $L=J+K$, and consider the algebra 
\[
B=E/L=E/(J+K) \cong (E/J)/((J+K)/J),
\] 
a quotient of $A=E/J$. Observe that, modulo the ideal $K$,  the generator $\Delta_{i,j}$ of $J$ reduces to $\Delta_{i,j}'=
a_{i,g}b_{i,g}+a_{j,g}b_{j,g}+ b_{i,1} a_{j,1}-a_{i,1} b_{j,1}$.

The ordering of the generators of the exterior algebra $E$ induces the degree-lexicographic order on the set of standard monomials in $E$. If the generators of $E$ are (generically) denoted $e_1,e_2,\dots,e_m$ and $S=(i_1,\dots,i_p)$ is an increasingly ordered subset of $\{1,\dots,m\}$, write $e_S=e_{i_1}\cdots e_{i_p}$ for the corresponding standard monomial in $E$. If $T=(j_1,\dots,j_q)$, then $e_S<e_T$ if $p<q$, or if $p=q$ and there exists $k$ with $1\le k\le p$ so that $i_r=j_r$ for $r<k$ and $i_k<j_k$. This order is multiplicative in the sense that if $e_S$ and $e_T$ are nontrivial standard monomials with $e_Se_T\neq 0$, then $e_Se_T$ is a standard monomial up to sign, and $1<e_S<\pm e_S e_T$. 

If $f=\sum c_S e_S$ is an element of $E$, the initial term $\mathrm{in}(f)$ of $f$ is the term $c_S e_S$ in this sum for which $e_S$ is the largest monomial among all $S$ for which $c_S\neq 0$. This monomial is the initial monomial of $f$, $e_S=\mathrm{inm}(f)$. For an ideal $F$ of $E$, the initial ideal $\mathrm{in}(F)$ of $F$ is the ideal generated by the initial terms $\mathrm{in}(f)$, $f\in F$. A set of elements $f_1,\dots,f_s \in F$ is a Gr\"obner basis of $F$ if 
$\mathrm{in}(f_1),\dots,\mathrm{in}(f_t)$ generate the ideal $\mathrm{in}(F)$.

For the specific exterior algebra $E$ generated by $a_{i,p},b_{i,p}$, $1\le i\le n$, $1\le p \le g$, above, order the generators as follows.
\[
\begin{aligned}
a_{1,1}>b_{1,1}>\dots>a_{1,g}>b_{1,g}>a_{2,1}>b_{2,1}&>\dots>a_{2,g}>b_{2,g}>\cdots\\
&\dots>a_{n,1}>b_{n,1}>\dots>a_{n,g}>b_{n,g}.
\end{aligned}
\]
\begin{prop} \label{prop:ggb}
The set
\[
\mathcal G=\left\{
\begin{matrix*}[l]a_{i,r}a_{i,s},\ a_{i,r}b_{i,s},\ b_{i,r}b_{i,s},&(1\le i\le n, 1\le r\neq s\le g),\\ 
a_{i,k}b_{i,k}-a_{i,g}b_{i,g},&(1\le i\le n, 1\le k\le g-1),\\
a_{i,1} b_{j,1}-b_{i,1} a_{j,1}-a_{i,g}b_{i,g}-a_{j,g}b_{j,g},&(1\le i<j\le n),\\
 a_{i,p}a_{j,q},\ a_{i,p}b_{j,q},\ b_{i,p}b_{j,q},&(1\le i \neq j\le n, 2\le p,q \le g),\\
 a_{i,1}a_{j,g}b_{j,g}+a_{i,g}b_{i,g}a_{j,1},&(1\le i<j\le n),\\
 b_{i,1}a_{j,g}b_{j,g}+a_{i,g}b_{i,g}b_{j,1},&(1\le i<j\le n)
\end{matrix*}\right\}
\]
is a Gr\"obner basis for the ideal $L=J+K$ in the exterior algebra $E$.
\end{prop}
\noindent Note that the elements of $\mathcal G$ are recorded with their initial terms first, and  that we have used $-\Delta_{i,j}'$ in place of $\Delta_{i,j}'$. Note also the presence of the cubic elements in~$\mathcal G$.

The proposition may be established using \cite[Corollary 1.5]{AHH}, by showing that all $S$- and $T$-polynomials involving elements of $\mathcal G$ reduce to zero with respect to $\mathcal G$.
This (lengthy) process may be inductively sped up, by successively considering the ideals $L^k$ and sets $\mathcal G^k$ involving the generators $a_{i,p}$, $b_{i,p}$ of $E$ with first index $i\ge k$, for $k=n-1,n-2,\dots$.

\begin{remark}
Gonz\'alez-Guti\'errez \cite{dc_GG} consider a further quotient of the algebra $A$. In the above notation, they work with the algebra $E/(L+K')$, where $K'$ is the ideal generated by the elements 
$(a_{i,1}-a_{1,1})(b_{j,1}-b_{1,1})$, $2\le i \neq j \le n$. This simplifies the zero divisor calculations carried out in \cite{dc_GG}, but complicates the above Gr\"obner basis considerations.
\end{remark}

\subsection{Zero divisors} \label{subsec:zd} 
The algebras $A=E/J$ and $B=E/L=E/(J+K)\cong A/((J+K)/J)$ we consider are quotients of the exterior algebra $E$ by ideals generated in degrees greater than or equal to two, so we identify the degree-one generators of all of these algebras and denote them by the same symbols $a_{i,p}$, $b_{i,p}$, $1\le i \le n$, $1\le p \le g$. Let $Z_A$ be the ideal in $A\otimes A$ generated by all degree-one zero divisors
\[
Z_A=\langle \bar{a}_{i,p}=1\otimes a_{i,p}-a_{i,p}\otimes 1,\ \bar{b}_{i,p}=1\otimes b_{i,p}-b_{i,p}\otimes 1, \quad (1\le i \le n, 1\le p \le g)\rangle.
\]
Similarly, denote the ideal generated by degree-one zero divisors in $B\otimes B$ by $Z_B$. 
To show that $\zcl(A) \ge 2n+2$, it suffices to show that $Z_A^{2n+2}\neq 0$. 
\begin{prop} \label{prop:prin}
The ideal $Z_A^{2n+2}$ is nonzero in $A\otimes A$.
\end{prop}
Let 
$\zeta=\bar{a}_{1,1}\bar{b}_{1,1}\bar{a}_{1,g}\bar{b}_{1,g}\bar{a}_{2,1}\bar{b}_{2,1}\cdots \bar{a}_{n,1} \bar{b}_{n,1}$. 
We assert that the image $\zeta_A$ of $\zeta$ is nonzero in $A\otimes A$. Clearly, this image is in $Z_A^{2n+2}$. For $n=1$, this is immediate, as $A=H^*(\Sigma_g)$ and $\zeta=\bar{a}_{1,1}\bar{b}_{1,1}\bar{a}_{1,g}\bar{b}_{1,g}=2a_{1,g}b_{1,g}\otimes a_{1,g}b_{1,g}=2\omega_1\otimes \omega_1$ in this instance. For $n\ge 2$, as the natural projection $A\otimes A \twoheadrightarrow B \otimes B$ takes the generators and powers of $Z_A$ to those of $Z_B$, it is enough to show that the image $\zeta_B$ of $\zeta$ in $B\otimes B$ is nonzero.

A calculation with the description of the ideal $L=J+K$ defining $B=E/L$ given in \S\ref{subsec:gb} reveals that, for $n\ge 2$, the image of $\zeta$ in $B\otimes B$ is given by
\[
\zeta_B=
2\sum_{k=0}^{n-1} (-1)^{n-k-1} \binom{n-1}{k} a_{1,g}b_{1,g}U_k \otimes a_{1,g}b_{1,g}V_{k},
\]
where $U_k=b_{2,1}\cdots b_{k+1,1}a_{k+2,1}\cdots a_{n,1}$ and $V_k=a_{2,1}\cdots a_{k+1,1}b_{k+2,1}\cdots b_{n,1}$. In particular, $U_0=V_{n-1}=a_{2,1}\cdots a_{n,1}$ and $V_0=U_{n-1}=b_{2,1}\cdots b_{n,1}$.

The Gr\"obner basis $\mathcal G$ recorded in Proposition \ref{prop:ggb} may be used to show that this element is nonzero in $B\otimes B$. For instance, the leading term of this element (in $E\otimes E$) is the tensor product $a_{1,g}b_{1.g}a_{2,1}\cdots a_{n,1}\otimes a_{1,g}b_{1.g}b_{2,1}\cdots b_{n,1}$ of two monomials neither of which reduce to zero with respect to $\mathcal G$. It follows that (the leading term of) $\zeta_B$ is nonzero in $B\otimes B$. Consequently, $\zeta_A\neq 0$, $Z_A^{2n+2}\neq 0$, and $\zcl(A) \ge 2n+2$ as was required.

\section{Orbit configuration spaces}

Let $\Gamma$ be a group and $X$ a $\Gamma$-space. The orbit configuration space $F_\Gamma(X,n)$ is the space of all ordered $n$-tuples of points in $X$ which lie in distinct $\Gamma$-orbits,
\[
F_\Gamma(X,n)=\{(x_1,\dots,x_n)\mid \Gamma\cdot x_i \cap \Gamma\cdot x_j=\emptyset\ \text{if}\ i\neq j\}.
\]
Orbit configuration spaces, introduced by Xicot\'encatl \cite{dc_X}, are natural generalizations of classical configuration spaces. If $\Gamma=\{1\}$ is trivial, $F_{\{1\}}(X,n)=F(X,n)$ is the classical configuration space.

\subsection{Generalized Fadell-Neuwirth theorem}
We will focus on the case where $X$ is a connected manifold without boundary of positive dimension, and $\Gamma$ is a finite group acting freely on $X$. Let ${\mathcal O}_n^{\Gamma}$ denote the union of $n$ distinct orbits, $\Gamma\cdot x_1,\dots,\Gamma\cdot x_n$, in $X$. The Fadell-Neuwirth theorem recorded in Theorem \ref{thm:FN} 
was generalized by Xicot\'encatl to orbit configuration spaces as follows.

\begin{theorem}[\cite{dc_X}] \label{thm:X} For $\ell \le n$, the projection onto the first $\ell$ coordinates, 
\[p_\Gamma \colon F_\Gamma(X,n)\to F_\Gamma(X,\ell),\] is a locally trivial bundle, with fiber 
$F_\Gamma(X\smallsetminus {\mathcal O}_n^\Gamma,n-\ell)$.
\end{theorem}
The proof of this result given in \cite{dc_X} is a modification of that of \cite{dc_FN} for classical configuration spaces. In the special case $\ell=n-1$ an alternative argument, which informs on the structure of these bundles, is given in \cite{dc_mono}.

Assume that the order of the finite group $\Gamma$ is $r$, and write $\Gamma=\{g_1,\dots,g_r\}$. Define a map from the orbit configuration space to the classical configuration space by sending an $n$-tuple of points in $X$ to their orbits. That is, define $f\colon F_\Gamma(X,n) \to F(X,rn)$ by \[
f(x_1,\dots,x_n) = (g_1x_1,\dots,g_rx_1,g_1x_2,\dots,g_rx_2,\dots\dots,g_1x_n,\dots,g_rx_n).
\]

\begin{theorem}[\cite{dc_mono}] \label{thm:mono}
The orbit configuration space bundle $p_\Gamma\colon F_\Gamma(X,n+1) \to F_\Gamma(X,n)$ is equivalent to the pullback of the classical configuration space bundle $p\colon F(X,rn+1) \to F(X,rn)$ under the map $f$.
\end{theorem}

Now specialize to the case where the finite cyclic group $\Gamma=\Z/r\Z$ acts freely on the manifold $X=\C^*$ by multiplication by the primitive $r$-th root of unity $\zeta=\exp(2\pi i/r)$. The associated orbit configuration space is 
\[
F_{\Gamma}(\C^*,n)=\{(x_1,\dots,x_n)\in(\C^*)^n \mid x_j \neq \zeta^k x_i,\ i\neq j,\ 1\le k\le r\},
\]
which may be realized as the complement in $\C^n$ of the hyperplane arrangement $\A_{r,n}$ consisting of the hyperplanes $H_i=\ker(x_i)$, $1\le i \le n$, and $H_{i,j}^k=\ker(x_i-\zeta^k x_j)$, $1\le i<j\le n$, $1\le k \le r$.
The arrangement $\A_{r,n}$ consists of the reflecting hyperplanes of the full monomial group,
the complex reflection group isomorphic to the wreath product of the symmetric group $\mathfrak{S}_n$ and $\Gamma=\Z/r\Z$. For instance, when $r=2$, this is the type B Coxeter group, and $\pi_1(F_{\Z/2\Z}(\C^*,n))$ is the type B pure braid group. Discussions of reflection arrangements, including the full monomial arrangements $\A_{r,n}$, may be found in references including \cite{OT,YuzHTC}.

\begin{theorem} $\tc(F_{\Z/r\Z}(\C^*,n))=2n$ \label{thm:orb}
\end{theorem}

This result may be established using techniques from the theory of hyperplane arrangements as discussed in Section \ref{sec:plane} and below, or by using the group theoretic methods presented in Section \ref{sec:groups}.

\subsection{Arrangements, II} \label{subsec:arr2}
Beginning with work of Arnol'd and Brieskorn (see \S\ref{subsec:cohomology}), 
the cohomology ring of the complement of a complex hyperplane arrangement is a well-studied object, facilitating analysis of the (higher) zero divisor cup length in this context. Let $\A=\{H_1,\dots,H_m\}$ be an arrangement of $m$ hyperplanes in $V=\C^\ell$. For convenience, we will assume that $\A$ is essential, and we will use cohomology with coefficients in $\C$.

The Orlik-Solomon theorem \cite{OS} shows that $H^*(M(\A);\C)$ is isomorphic to the Orlik-Solomon algebra $A(\A)$, the quotient of the exterior algebra $E(\A)$ generated by one-dimensional classes $e_j$, $1\le i \le m$, by a homogeneous ideal $I(\A)$. Detailed expositions may be found in \cite{OT,YuzOS}. Let $[m]=\{1,\dots,m\}$, refer to the hyperplanes of $\A$ by their subscripts, and order them accordingly. Given $S\subset [m]$, denote the flat $\bigcap _{i\in S} H_i$ by $\cap S$. If  $\cap S\neq \emptyset$, call $S$ independent if the codimension of $\cap S$ in $V$ is equal to $|S|$, and dependent if $\mathrm{codim}(\cap S)<|S|$. If $S=(i_1,i_2,\dots,i_p)$ is an increasingly ordered subset of $[m]$, recall that $e_S=e_{i_1}e_{i_2}\cdots e_{i_p}$ denotes the corresponding standard monomial in the exterior algebra. Define $\partial e_S=\sum_{k=1}^p (-1)^{k-1} e_{S\setminus\{i_k\}}$. The Orlik-Solomon ideal is generated by
\[
\{\partial e_S \mid S\ \text{is dependent}\}\bigcup \{e_S \mid \cap S=\emptyset\}.
\]

A circuit is a minimally dependent subset $T\subseteq[m]$, that is, $T$ is dependent, but every nontrivial subset of $T$ is independent. If $S$ is dependent and $T\subset S$ a circuit, then $e_S=\pm e_T e_{S\setminus T}$, and $\partial e_S=\pm \partial e_T\cdot e_{S\setminus T} \pm e_T\cdot\partial e_{S\setminus T}$. Also, note that if $\cap T=\emptyset$, then $\cap S=\emptyset$ for any $S$ containing $T$. 
In light of these two observations, the generating set for the Orlik-Solomon ideal given above may be reduced as follows: The ideal $I(\A)$ is generated by
\begin{equation} \label{eq:OSgb}
\mathcal G=\{\partial e_T \mid T\ \text{is a circuit}\} \bigcup \{e_S \mid S\ \text{is minimal such that}\ \cap S=\emptyset\}.
\end{equation}

The ordering of the hyperplanes of $\A$ induces the degree-lexicographic order on the set of standard monomials $e_S$ in the exterior algebra $E(\A)$. Call a subset $S$ of $[m]$ a broken circuit if there exists $k\in[m]$ so that $k<i$ for all $i\in S$ and $(k,S)$ is a circuit. Broken circuits correspond to the initial monomials of the elements $\partial e_T$ appearing in \eqref{eq:OSgb}. In \cite[Theorem 2.8]{YuzOS}, Yuzvinsky shows that the initial monomials  of elements of $\mathcal G$ generate $\mathrm{in}(I(\A))$, the ideal generated by the initial terms of elements of $I(\A)$, whence $\mathcal G$ is a Gr\"obner basis for the Orlik-Solomon ideal $I(\A)$, see \cite{AHH}. This yields a basis for the quotient $A(\A)=E(\A)/I(\A)$, the Orlik-Solomon algebra of the arrangement $\A$.

For $S\subset[m]$, let $a_S$ denote the image of the standard monomial $e_S$ in the Orlik-Solomon algebra $A(\A)$. The \textbf{nbc} basis for $A(\A)$ consists of all elements $a_S$ corresponding to subsets $S$ of $[m]$ which contain \textbf{n}o \textbf{b}roken \textbf{c}ircuits \cite{OT,YuzOS}. This basis has been used to study the (higher) zero divisor cup length of $A(\A)$ by a number of authors, including \cite{dc_FY,YuzHTC}. Following these references, we restrict our attention to a central arrangement $\A$. Recall that we assume $\A$ is essential in $\C^\ell$. Consequently, a maximal independent set (resp., \textbf{nbc} set) has cardinality $\ell$.

Let $\Pi=(B,C)$ be an ordered pair of disjoint subsets of $[m]$, and let $\bar{\Pi}=B\cup C$. The pair $\Pi$ is said to be basic if $B$ and $C$ are \textbf{nbc} sets for some linear order on $\bar{\Pi}$ and $B$ is a maximal independent set, $|B|=\ell$. The central arrangement $\A$ is said to be large if there is a basic pair $\Pi=(B,C)$ with $|C|=\ell-1$. In \cite{YuzHTC}, Yuzvinsky uses basic pairs to find lower bounds on the (higher) zero divisor cup length of the Orlik-Solomon algebra, and proves the following.

\begin{theorem}[\cite{YuzHTC}] \label{thm:large}
Let $\A$ be an essential central arrangement in $\C^\ell$ with complement $M(\A)$, and let $s\ge 2$ be a positive integer. If $(B,C)$ is a basic pair, then $\tc_s(M(\A)) > (s-1)\ell+|C|$. If $\A$ is large, then $\tc_s(M(\A)) = s\ell$.
\end{theorem}

The arrangements $\A_{r,n}$ in $\C^n$ associated to the full monomial groups and arising in the context of cyclic group orbit configuration spaces are large. Recall that $\A_{r,n}$ has hyperplanes $H_i=\ker(x_i)$ and $H_{i,j}^k=\ker(x_i-\zeta^k x_j)$, where $\zeta=\exp(2\pi i/r)$, and take $B=\{H_1,\dots,H_n\}$ and $C=\{H_{1,2},H_{1,3},\dots,H_{1,n}\}$. Thus, Theorem \ref{thm:orb} follows from Theorem \ref{thm:large}. More generally, Yuzvinsky establishes an analogous result for the reflection arrangement associated to 
any irreducible complex reflection group.

A complex reflection in $V=\C^\ell$ is a finite order linear transformation $\tau\colon V \to V$ whose fixed point set is a hyperplane $H_\tau$. A reflection group is a finite subgroup of $\mathrm{GL}(V)$ that is generated by reflections. A reflection group is irreducible if its tautological representation in 
$\mathrm{GL}(V)$ is irreducible. The reflection arrangement $\A_W$ associated to the reflection group $W$ is the set of hyperplanes $\{H_\tau \mid \tau\ \text{a reflection in}\ W\}$.
\begin{theorem}[\cite{YuzHTC}] 
Let $W$ be an irreducible reflection group of rank $\ell$, and let $s\ge 2$ be a positive integer. If $\A_W$ is the associated reflection arrangement, then $\tc_s(M(\A_W))=s\ell$.
\end{theorem}

\begin{remark} If $\varGamma$ is a simple graph with vertices $\{1,\dots,n\}$, the associated graphic arrangement $\A_\varGamma$ consists of the hyperplanes $\ker(x_i-x_j)$ in $\C^n$ corresponding to the edges $\{i,j\}$ of $\varGamma$. For instance, if $\varGamma=K_n$ is the complete graph, then $\A_\varGamma=\A_n$ is the  braid arrangement introduced in Example \ref{ex:braid arr}. In \cite{dc_NF}, Fieldsteel uses Yuzvinsky's result stated in Theorem \ref{thm:large} to find conditions on the graph $\varGamma$, related to the arboricity, which insure that the (higher) topological complexity of the complement $M(\A_\varGamma)$ of a graphic arrangement  is as large as possible.
\end{remark}

\section{Some discrete groups}
\label{sec:groups}

Let X be an aspherical space, that is, a space whose higher homotopy groups vanish: $\pi_i(X) = 0$ for $i \ge 2$. Farber \cite{dc_Fa05} poses the problem of computing the topological complexity of such a space in terms of algebraic properties of the fundamental group $G = \pi_1(X)$. In other words, given a discrete group $G$, define the topological complexity of $G$ to be $\tc(G) := \tc(K(G,1))$, the topological complexity of an Eilenberg-Mac\,Lane space of type $K (G, 1)$, and express $\tc(G)$ in terms of invariants such as the cohomological or geometric dimension of $G$ if possible.

\begin{example} Associated to a simple graph $\varGamma$ on $n$ vertices is a right-angled Artin group $G_\varGamma$ with generators corresponding to the vertices of $\varGamma$, and commutator relators corresponding to the edges. For instance, if $\varGamma=K_n$ is the complete graph, then $G_\varGamma=\Z^n$ is free abelian, while if $\varGamma$ has no edges, then $G_\varGamma=F_n$ is free. For any right-angled Artin group, one has $\tc(G_\varGamma)=z(\varGamma)+1$, where $z(\varGamma)$ is the maximal number of vertices of $\varGamma$ covered by two (disjoint) cliques in $\varGamma$, see \cite{cp,ggy,dc_GLO}.
\end{example}

Many of the configuration spaces discussed previously are $K(G,1)$-spaces, for surface pure braid groups, for pure braid groups associated to reflection groups, etc. For example, $\pi_1(F(\C,n))=P_n$ is the Artin pure braid group. From the homotopy exact sequence of the Fadell-Neuwirth bundle $F(\C,m) \to F(\C,m-1)$, with fiber $\C\smallsetminus\{m-1\ \text{points}\}$ and cross section, we see (inductively) that $F(\C,n)$ is a $K(P_n,1)$-space, and obtain a split, short exact sequence $1\to F_{n-1} \to P_n \to P_{n-1} \to 1$, where $F_k$ is the free group on $k$ generators. Thus,
\[
P_n = F_{n-1} \rtimes P_{n-1} = 
F_{n-1} \rtimes (F_{n-2} \rtimes P_{n-2})=
\dots = F_{n-1} \rtimes ( \dots \rtimes (F_3\rtimes (F_2\rtimes F_1)))
\]
is an iterated semidirect product of free groups. 

The iterated semidirect product structure of $P_n$ is apparent in the classical presentation of this group. The pure braid group $P_n$ has generators $A_{i,j}$, $1\le i<j \le n$, and relations
\begin{equation}\label{eq:Prels}
A_{r,s}^{-1}A_{i,j}^{}A_{r,s}^{}=\begin{cases}
A_{i,j}&\text{if $r<s<i<j$,}\\
A_{r,j}^{}A_{s,j}^{}A_{r,j}^{}A_{s,j}^{-1}A_{r,j}^{-1}&\text{if $r=i<s<j$,}\\
A_{r,j}^{}A_{s,j}^{}A_{r,j}^{-1}&\text{if $r<i=s<j$,}\\
[A_{r,j},A_{s,j}]A_{i,j}[A_{r,j},A_{s,j}]^{-1}&\text{if $r<i<s<j$,}\\
A_{i,j}&\text{if $i<r<s<j$,}
\end{cases}
\end{equation}
where $[u,v]=uvu^{-1}v^{-1}$ denotes the commutator, see, for instance Birman \cite{dc_Bi}. 
Observe that, for $s<j$ as in the relations above, the action of $F_{s-1}=\langle A_{1,s},\dots,A_{s-1,s}\rangle$ on 
$F_{j-1}=\langle A_{1,j},\dots,A_{j-1,j}\rangle$ (via the Artin representation) is by conjugation. 
It follows that the induced action of $P_{n-1}$ on $H_*(F_{n-1},;\Z)$ is trivial.

\subsection{Almost-direct products of free groups} \label{subsec:adfg}
An \emph{almost-direct product of free groups} is an iterated semidirect product $G=F_{d_n} \rtimes \cdots \rtimes F_{d_1}$ of finitely generated free groups for which $F_{d_i}$ acts trivially on $H_*(F_{d_j};\Z)$ for $i<j$. Thus, $P_n$ is an almost-direct product of free groups. The fundamental groups of the orbit configuration spaces $F_{\Z/r\Z}(\C^*,n)$ considered in the previous section provide another family of examples.

Let $\Gamma=\Z/r\Z$, and $\zeta=\exp(2\pi i/r)$. The pure braid group $P_{r,n}=\pi_1(F_{\Gamma}(\C^*,n))$ associated to the full monomial group $G(r,n)$ may also be realized as an almost-direct product of free groups. From Theorems \ref{thm:X} and \ref{thm:mono}, the map $F_{\Gamma}(\C^*,n) \to F_{\Gamma}(\C^*,n-1)$ defined by forgetting the last coordinate is a bundle, with fiber $\C^* \smallsetminus \{n-1\ \text{orbits}\}=\C\smallsetminus\{r(n-1)+1\ \text{points}\}$. A minor modification of these results is useful in revealing the almost-direct product structure of $P_{r,n}$. Given a configuration of $m$ distinct ordered points $(x_1,\dots,x_m)$ in $\C^*$, one obtains a configuration of $m+1$ distinct ordered points $(0,x_1,\dots,x_m)$ in $\C$, yielding a homotopy equivalence $F(\C^*,m)\simeq F(\C,m+1)$. Using this observation, together with Theorem \ref{thm:mono}, one can check that the bundle $F_{\Gamma}(\C^*,n) \to F_{\Gamma}(\C^*,n-1)$ may be realized as the pullback of the classical configuration space bundle $F(\C,N+1)\to F(\C,N)$ where $N=r(n-1)+1$, under the map $g\colon F_\Gamma(\C^*,n-1) \to F(\C,N)$ given by
\[
g(x_1,\dots,x_{n-1})=(0,\zeta x_1,\dots,\zeta^r x_1,\zeta x_2,\dots,\zeta^r x_2,\dots\dots,\zeta x_{n-1},\dots,\zeta^r x_{n-1}).
\]
It follows that the orbit configuration space bundle $F_{\Gamma}(\C^*,n) \to F_{\Gamma}(\C^*,n-1)$ admits a section, and the fundamental group $P_{r,n-1}$ of the base acts trivially on the homology of the fiber.  Hence, an inductive argument reveals that $P_{r,n}$ is an almost-direct product of free groups.

Under natural assumptions on the ranks of the constituent free groups, the topological complexity of an almost-direct product of free groups was determined by Cohen.

\begin{theorem}[\cite{dc_ad}] \label{thm:ad}
 If $G=F_{d_n} \rtimes \cdots \rtimes F_{d_1}$ is an almost-direct product of free groups with $d_j \ge 2$ for each $j$, and $m$ is a nonnegative integer, then
 \[
 \tc(G\times \Z^m)=2n+m+1.
 \]
\end{theorem}
For an arbitrary iterated semidirect product of free groups $G= \rtimes_{j=1}^n F_{d_j}$ of cohomological dimension $n$, a $K(G,1)$-complex of dimension $n$ is constructed in \cite{dc_CS}. Thus, for such groups, the dimensional upper bound on topological complexity may be stated in terms of the cohomological dimension as $\tc(G) \le 2 \cd(G)+1$.

The integral homology $H_*(G;\Z)$ is torsion-free and the Poincar\'e polynomial is given by $P(G,t) = \sum_{k=0}^n b_k(G) \cdot t^k = \prod_{j=1}^n (1+d_j t)$, where $b_k(G)$ is the $k$-th Betti number of $G$, see \cite{dc_FR}.   A minimal, free $\mathbb Z{G}$-resolution of $\mathbb Z$, which we denote by $C_\bullet(G) \xrightarrow{\ \epsilon\ } \mathbb Z$, is constructed in \cite{dc_CS}.  

Let $N=b_1(G)=d_1+d_2+\dots+d_n$.  The abelianization map $\mathfrak a\colon G \to \mathbb Z^N$ induces a chain map $\mathfrak a_\bullet\colon 
C_\bullet \to K_\bullet$, where $C_\bullet=C_\bullet(G) \otimes_{\mathbb Z{G}} \mathbb Z \mathbb Z^N$ and $K_\bullet \to \mathbb Z$ is the standard  $\mathbb Z \mathbb Z^N$-resolution of $\mathbb Z$.  The induced map $\mathfrak a^*\colon H^*(\Z^N;\Z) \to H^*(G;\Z)$ in integral cohomology is 
surjective, and is an isomorphism $\mathfrak a^*\colon H^1(\Z^N;\Z) \xrightarrow{\ \sim\ } H^1(G;\Z)$ in dimension one, see \cite[Theorem 2.1]{dc_ad}.

Let $J$ be the ideal in the exterior algebra $H^*(\Z^N;\Z)$ generated by the elements of the kernel of the surjection $\mathfrak a^2 \colon H^2(\Z^N;\Z) \to H^2(G;\Z)$, $J=(\ker(\mathfrak a^2))$. An explicit Gr\"obner basis for $J$ is exhibited in \cite[\S3]{dc_ad} (in the degree-lexicographic order on a standard basis for the exterior algebra), and this is used to shown that the integral cohomology ring of $G$ is given by $H^*(G;\Z) \cong H^*(\Z^N;\Z)/J$. 

Passing to field coefficients, $H^*(-)=H^*(-;\C)$, if $G=\rtimes_{j=1}^nF_{d_j}$ is an almost-direct product of free groups with $d_j \ge 2$ for each $j$, one can exhibit pairs of generators of $x_i,y_j\in H^1(G)$ corresponding to distinct generators of the free groups $F_{d_j}$, $1\le j \le n$.  As shown in \cite[Theorem 4.2]{dc_ad}, this yields $2n$ zero-divisors $\bar{x}_j,\bar{y}_j$ in $H^1(G)\otimes H^1(G)$ with nonzero product. These considerations yield $\tc(G)=2n+1$ for $G$ as in the statement of the theorem. The general case $\tc(G\times\Z^m)=2n+m+1$ may be obtained from this, the product inequality, and a straightforward analysis of the zero-divisor cup length of $H^*(G\times \Z^m)$.

\subsection{Fiber-type arrangements}
The Artin pure braid group $P_n$ associated to the symmetric group, and, more generally, the pure braid groups $P_{r,n}$ associated to  the full monomial groups 
may be realized as the fundamental groups of the hyperplane arrangements defined by the polynomials 
\[
Q(\A_{r,n})=x_1\cdots x_n \cdot \prod_{1\le i<j \le n}(x_i^r-x_j^r), \ r\ge 1.
\]
Notice that the arrangement $\A_{1,n}$ here corresponds to the arrangement $\hat\A_{n+1}$ from Example \ref{ex:braid arr}, so that the fundamental group $\pi_1(M(\A_{1,n}))=P_{n+1}$ is the Artin pure braid group on $n+1$ strands.

The arrangements $\A_{r,n}$ are examples of (essential) fiber-type arrangements. An arrangement $\A$ in $\C^{\ell+1}$ is said to be strictly linearly fibered if there is a choice of coordinates $(\mathbf{x},z)=(x_1,\dots,x_\ell,z)$ on $\C^{\ell+1}$ so that the restriction, $\varrho$, of the projection $\C^{\ell+1}\to\C^\ell$, $(\mathbf{x},z)\mapsto \mathbf{x}$, to the complement $M(\A)$ is a fiber bundle projection, with base $\varrho(M(\A))=M(\mathcal B)$, the complement of an arrangement $\mathcal B$ in $\C^\ell$, and fiber the complement of finitely many points in $\C$. We say $\A$ is strictly linearly fibered over $\mathcal B$. Fiber-type arrangements are then defined inductively as follows: An arrangement $\A=\A_1$ of finitely many points in $\C^1$ is fiber-type. An arrangement $A=\A_\ell$ of hyperplanes in $\C^\ell$ is fiber-type if $\A$ is strictly linearly fibered over a fiber-type arrangement $\A_{\ell-1}$ in $\C^{\ell-1}$.

The complement of a fiber-type arrangement sits atop a tower of fiber bundles
\[
M(\A_\ell) \xrightarrow{\varrho_\ell} M(\A_{\ell-1}) \xrightarrow{\varrho_{\ell-1}} \cdots \xrightarrow{\varrho_2} M(\A_1)=\C\smallsetminus\{d_1\ \text{points}\},
\]
where the fiber of $\varrho_k$ is homeomorphic to the complement of $d_k$ points in $\C$. Repeated application of the homotopy exact sequence of a bundle shows that $M(\A_\ell)$ is a $K(\pi,1)$-space, where $\pi=\pi_1(M(\A_\ell))$. The integers $\{d_1,\dots,d_\ell\}$ are called the exponents of the fiber-type arrangement $\A_\ell$.

Suppose $\A$ is strictly linearly fibered over $\mathcal B$, and $|\A|=m+n$, where $|\mathcal B|=m$. From the definition, a defining polynomial for $\A$ factors as $Q(\A)=Q(\mathcal B) \cdot \phi(\mathbf{x},z)$, where $\phi(\mathbf{x},z)=\prod_{j=1}^n(z-g_j(\mathbf{x})$ is a product of $n$ linear functions. 
Since $\phi(\mathbf{x},z)$ has distinct roots for any $\mathbf{x}\in M(\mathcal B)$, the map
\[
g\colon M(\mathcal B) \to \C^n,\quad g(\mathbf{x})=\bigl(g_1(\mathbf{x}),g_2(\mathbf{x}),\dots,g_n(\mathbf{x})\bigr),
\]
takes values in the configuration space $F(\C,n)$.

\begin{theorem}[\cite{dc_mono}] \label{thm:ft}
Let $\mathcal B$ be an arrangement of $m$ hyperplanes, and let $\A$ be an arrangement of $m+n$ hyperplanes which is strictly linearly fibered over $\mathcal B$. Then the bundle $\varrho\colon M(\A) \to M(\mathcal B)$ is equivalent to the pullback of the classical configuration space bundle $p\colon F(\C,n+1) \to F(\C,n)$ under the map $g$.
\end{theorem}
From this result, it follows that the bundle $\varrho\colon M(\A) \to M(\mathcal B)$ admits a section, that the structure group of this bundle is the pure braid group $P_n$, and that the fundamental group of the base acts by conjugation (in fact, by a pure braid action) on the fundamental group of the fiber. 

If $\A=\A_\ell$ is a fiber-type arrangement with exponents $\{d_1,\dots,d_\ell\}$, repeated application of this theorem and these consequences reveals that 
\[
\pi_1(M(\A_\ell))=F_{d_\ell} \rtimes \cdots \rtimes F_{d_1}
\]
is an almost-direct product of free groups. Theorem \ref{thm:ad} yields the following.
\begin{cor} Let $\A_\ell$ be a fiber-type arrangement with exponents $\{d_1,\dots,d_\ell\}$ and let $G(\A_\ell)=\pi_1(M(\A_\ell))$ be the fundamental group of the complement of $\A_\ell$. If the exponents of $\A_\ell$ are all at least two, $d_j \ge 2$, $1\le j \le \ell$, and $m$ is a nonnegative integer, then
\[
\tc(G(\A_\ell) \times \Z^m) 
=2\ell+m+1.
\]
\end{cor}

Theorem \ref{thm:orb} may be obtained from this result as follows. From Proposition \ref{prop:central/affine} and Remark \ref{rem:decone}, the complement of the essential, central arrangement $\A_{r,n}$ may be realized as $M(\A_{r,n}) \cong M(\mathrm{d}\A_{r,n}) \times \C^*$. It is readily checked that the decone $\mathrm{d}\A_{r,n}$ (with respect to $\{x_1=0\}$) in $\C^{n-1}$ is fiber-type, with exponents $r+1,\, 2r+1,\dots,\, r(n-1)+1$ all at least two. So we have
\[
\begin{aligned}
\tc(F_{\Z_r}(\C^*,n))&=\tc(M(\A_{r,n}))=\tc(P_{r,n})=\tc(G(\A_{r,n}))\\
&=\tc(M(\mathrm{d}\A_{r,n})\times \C^*)=\tc(G(\mathrm{d}\A_{r,n})\times \Z)\\
&=2(n-1)+1+1=2n.
\end{aligned}
\]

\subsection{Subgroup conditions}
Several of the results on the topological complexity of discrete groups mentioned above may also be obtained by other (group-theoretic) means, as shown by Grant-Lupton-Oprea.
\begin{theorem}[\cite{dc_GLO}]\label{thm:glo}
Let $G$ be a discrete group. If $H$ and $K$ are subgroups of $G$ which satisfy $gHg^{-1} \cap K=\{1\}$ for every $g \in G$, then $\tc(G) \ge \cd(H\times K)+1$.
\end{theorem}

This may be used to recover the topological complexity of the pure braid group, $\tc(P_n)=\tc(F(\C,n))=2n-2$. 
The pure braid group 
$P_n$ has a free abelian subgroup $H\cong \Z^{n-1}$, generated in terms of the standard generators $A_{i,j}$ of $P_n$ by $A_{j,j+1}A_{j,j+2} \cdots A_{j,n}$, $1\le j \le n-1$, see, for instance, Birman \cite{dc_Bi}.  Let $K< P_n$ be the image of the right splitting in the 
exact sequence
$1\to F_{n-1}\to P_n \to P_{n-1}\to 1$. The subgroup $K$ consists of pure braids with trivial last strand, and is generated by $A_{i,j}$ with $j<n$. 
In \cite{dc_GLO}, it is shown geometrically that 
$gHg^{-1} \cap K=\{1\}$ $\forall\,g\in P_n$. This may also be established algebraically using the pure braid relations recorded in \eqref{eq:Prels} above. 
Consequently, $\tc(P_n) \ge \cd(H \times K) +1= (n-1) + (n-2)+1=2n-2$.

Theorem \ref{thm:glo} may additionally be used to recover the topological complexity of the pure monomial braid group $P_{r,n}$, recorded in Theorem \ref{thm:orb} in terms of that of the orbit configuration space $F_{\Z_r}(\C^*,n)$, a $K(P_{r,n},1)$-space. 
From Theorem \ref{thm:mono}  and the discussion in \S\ref{subsec:adfg}, the group $P_{r,n}$ may be expressed as an almost-direct product 
$P_{r,n} \cong F_N \rtimes P_{r,n-1}$, where $N=r(n-1)+1$. Using this, and 
the presentation of $P_{r,n}$ given in \cite[Theorem 2.2.4]{dc_mono}, one can exhibit a rank $n$ free abelian subgroup $H<P_{r,n}$ and a cohomological dimension $n-1$ subgroup $K<P_{r,n}$ satisfying $gHg^{-1}\cap K=\{1\}$ $\forall g\in P_{r,n}$. Consequently, $\tc(P_{r,n}) \ge \cd(H\times K)+1=n+(n-1)+1=2n$. 
We anticipate that this result may also be used to recover the topological complexity of other almost-direct products of free groups, such as the fundamental groups of complements of other fiber-type hyperplane arrangements.

In \cite{dc_GLO}, Grant-Lupton-Oprea also use Theorem \ref{thm:glo} to find the topological complexity of right-angled Artin groups, and strikingly, to show that $\tc(\mathcal{H})=5$ for Higman's acyclic group~$\mathcal{H}$.

\begin{remark} \label{rem:GR-M}
The pure braid group $P_n$ is the kernel of the homomorphism $\pi\colon B_n \to\mathfrak{S}_n$ from the full braid group to the symmetric group sending a braid to its induced permutation. If $S \le \mathfrak{S}_n$ is a subgroup, the preimage $B_n^S:=\pi^{-1}(S)$ is a subgroup of $B_n$ containing $P_n$. In \cite{GR-M}, Grant and Recio-Mitter study the (higher) topological complexity of subgroups of $B_n$ arising in this way. For a subgroup $S\le\mathfrak{S}_{n-k}\times\mathfrak{S}_k\le\mathfrak{S}_n$, they use Theorem \ref{thm:glo} to show that $\tc(B_n^S) \ge 2n-k$. In particular, this lower bound, together with the upper bound $\tc(G) \le \cd((G\times G)/\mathcal{Z}(G))$ of Grant \cite{Grant} for a torsion-free discrete group $G$ with center $\mathcal{Z}(G)$ embedded in $G\times G$ via the diagonal homomorphism, yields $\tc(B_n^S)=2n-2$ when $S \le \mathfrak{S}_{n-k} \times \{1\}^2$.
\end{remark}

\section{Sins of omission} \label{sec:sins}
We close with some brief remarks regarding two directions not discussed in the previous sections, including comments indicating the reasons for these omissions.

\subsection{Graph configuration spaces} \label{subsec:graph}
Investigation of the collision-free motion of automated guided vehicles on, for instance, a network of wires leads one to study the configuration spaces of distinct points on graphs (see \cite{GK}), and the topological complexity of these spaces. Let $\varGamma$ be a finite, connected graph. In sharp contrast to the behavior  discussed in the previous sections, the topological complexity of the ordered configuration space $F(\varGamma,n)$ of a graph is, in certain instances, independent of the number of particles $n$.  
We state results of Farber and Scheirer, which apply in the case where $\varGamma$ is a tree, along these lines below.

An essential vertex of $\varGamma$ is a vertex which has at least three incident edges. Let $m(\varGamma)$ denote the number of essential vertices of $\varGamma$. As noted in \cite[Lemma 10.1]{FarberGraphs2}, the configuration space $F(\varGamma,n)$ is connected if $\varGamma$ has at least one essential vertex. This is the case if $\varGamma$ is not homeomorphic to the closed interval $[0,1]$ (and $n\ge 2$) or to the circle $S^1$ (and $n\ge 3$).
If $\varGamma\not\cong S^1$, a result of Ghrist \cite{Ghrist} shows that the configuration space $F(\varGamma,n)$ has the homotopy type of a cell complex of dimension at most $m(\varGamma)$. Thus the dimension inequality recorded in the Introduction insures that $\tc(F(\varGamma,n)) \le 2 m(\varGamma)+1$ for any such graph, see Farber \cite{FarberGraphs1}.

\begin{theorem}[\cite{FarberGraphs2}] \label{thm:Fgraph}
Let $\varGamma$ be a tree not homeomorphic to the closed interval $[0,1]$, and let $n$ be an integer satisfying $n \ge 2 m(\varGamma)$. If $n=2$, assume in addition that $\varGamma$ is not homeomorphic to the letter ${\sf{Y}}$. Then, $\tc(F(\varGamma,n))=2m(\varGamma)+1$.
\end{theorem}

An arc in $\varGamma$ is a subspace homeomorphic to the interval $[0,1]$. Let $\mathcal{V}$ be a subset of the vertex set of $\varGamma$. Loosely speaking, a collection of oriented arcs is allowable for $\mathcal{V}$ if no vertex in $\mathcal{V}$ is an endpoint of any arc, and at each vertex in $\mathcal{V}$ there is at least one direction in which the orientations of the arcs do not ``cancel out.'' Refer to Scheirer \cite[Definition 2.3]{Scheirer} for a precise formulation.

\begin{theorem}[\cite{Scheirer}]  \label{thm:Scheirer}
Let $\varGamma$ be a tree with $m(\varGamma)\ge 1$, and let $k$ be the smallest integer for which there is a collection of $k$ oriented arcs in $\varGamma$ which is allowable for the collection of all vertices in $\varGamma$ which have exactly three incident edges. Let $k=0$ if there are no such vertices. If $n\ge 2m(\varGamma)+k$, then $\tc(F(\varGamma,n))=2m(\varGamma)+1$.
\end{theorem}

\subsection{Unordered configuration spaces} \label{subsec:unordered}
Given a space $X$, the symmetric group $\mathfrak{S}_n$ acts on $X^n$ by permuting coordinates. This restricts to a free action on the configuration space $F(X,n)$ of $n$ distinct ordered points in $X$. The orbit space $F(X,n)/\mathfrak{S}_n$ is the configuration space of $n$ distinct unordered points in $X$, with fundamental group the full braid group of $X$. These unordered configuration spaces are sometimes $K(G,1)$-spaces. For example, $\pi_1(F(\C,n)/\mathfrak{S}_n)=B_n$ is the Artin full braid group, and $F(\C,n)/\mathfrak{S}_n$ is an aspherical space, since its covering space $F(\C,n)$ is. On the other hand, the full braid group of the sphere $S^2$ has torsion (see \cite{dc_Bi}),  so does not have a finite $K(G,1)$.

More generally, if $W$ is a complex reflection group, with associated reflection arrangement $\A_W$, then the complement $M(\A_W)$ is an aspherical space. This fact has a lengthy history, relevant references include \cite{dc_Ar,dc_Br,Deligne,OT,Bessis}. As above, $W$ acts freely on $M(\A_W)$, and the orbit space $M(\A_W)/W$ is also aspherical. The fundamental groups $P_W=\pi_1(M(\A_W))$ and $B_W=\pi_1(M(\A_W)/W)$ are the pure and full braid groups for $W$.

The reader has likely noted that, prior to this point, the topological complexity of unordered configuration spaces, orbit spaces of reflection groups, and the associated full braid groups has not been mentioned. 
The reason for this is quite simple. To the best of our knowledge, very little is known concerning the topological complexity in any of these contexts. While group theoretic aspects of full braid groups and cohomological aspects of these groups, configuration spaces, and orbit spaces are well studied, as noted in \cite{GR-M},  the lower bounds provided by the zero divisor cup length and by the subgroup conditions of Theorem \ref{thm:glo} appear to be insufficient to determine the topological complexity. This is the case, in particular, for the Artin full braid group $B_n$, that is, for the unordered configuration space $F(\C,n)/\mathfrak{S}_n$.

One notable exception is the unordered configuration space $F(\varGamma,n)/\mathfrak{S}_n$ of a tree $\varGamma$. Under the assumptions of Theorem \ref{thm:Scheirer}, Scheirer \cite{Scheirer} shows that the topological complexity of the unordered configuration space is equal to that of the ordered configuration space, 
$\tc(F(\varGamma,n)/\mathfrak{S}_n) = \tc(F(\varGamma,n))=2m(\varGamma)+1$.

\section*{Acknowledgements}
This survey expands on lectures given at the Mathematisches Forschungsinstitut Oberwolfach Mini-Workshop \emph{Topological Complexity and Related Topics} in the Spring of 2016. We thank the organizers M. Grant, G. Lupton, and L. Vandembroucq for an enjoyable, interesting, and stimulating workshop. 
We also thank the MFO for its support and hospitality, and for providing a productive mathematical envirnoment. The algebraic geometry and commutative algebra system \emph{Macaulay2} \cite{M2} and the package \emph{Mathematica} (from Wolfram Research) were useful in the preparation of the manuscript. Finally, we thank the referee for pertinent remarks and suggestions, and we thank S. He, J. Kona, J. Strummer, M. Jones, P. Simonon, and T. Headon for inspiration.

\newcommand{\arxiv}[1]{{\texttt{\href{http://arxiv.org/abs/#1}{{arXiv:#1}}}}}

\newcommand{\MRh}[1]{\href{http://www.ams.org/mathscinet-getitem?mr=#1}{MR#1}}

\bibliographystyle{amsplain}

\begin{thebibliography}{99}
\bibitem{AHH} A. Aramova, J. Herzog, T. Hibi, \emph{Gotzmann theorems for exterior algebras and
combinatorics}, J. Algebra \textbf{191} (1997), 174--211;  \MRh{1444495}.

\bibitem{dc_Ar} V. Arnol'd, \emph{The cohomology ring of the group of dyed braids}, Math. Notes \textbf{5} (1969), 138--140;\\ \MRh{0242196}.

\bibitem{Bessis} D. Bessis, \emph{Finite complex reflection arrangements are $K(\pi,1)$}, Ann. of Math. (2) \textbf{181} (2015), 809--904; \MRh{3296817}. 

\bibitem{dc_Bi} J. Birman, \emph{Braids, Links and Mapping Class Groups}, Ann. of Math. Studies No. 82,  Princeton Univ. Press, Princeton, NJ 1974; \MRh{0375281}.

\bibitem{dc_Br} E. Brieskorn, \emph{Sur les groupes de tresses}, in: \emph{S\'eminaire Bourbaki, 1971/72}, pp. 21--44, Lecture Notes in Math., vol. 317, Springer, Berlin, 1973; \MRh{0422674}.

\bibitem{dc_mono} D. Cohen, \emph{Monodromy of fiber-type arrangements and orbit configuration 
spaces}, Forum Math. \textbf{13} (2001), 505--530; \MRh{1830245}.

\bibitem{dc_ad} D. Cohen, \emph{Cohomology rings of almost-direct products of free groups},
Compos. Math. \textbf{146} (2010), 465--479; \MRh{2601635}.

\bibitem{dc_CF} D. Cohen, M. Farber, \emph{Topological complexity of collision-free motion planning on surfaces}, Compos. Math. \textbf{147} (2011), 649--660; \MRh{2776616}.

\bibitem{cp} D. Cohen, G. Pruidze, \emph{Motion planning in tori}, Bull. Lond. Math. Soc. \textbf{40} (2008), 249--262;\\ \MRh{2414784}.

\bibitem{dc_CS} D. ~Cohen, A. ~Suciu, \emph{ Homology of iterated semidirect
products of free groups}, J. Pure Appl. Algebra \textbf{126} (1998),
87--120; \MRh{1600518}.

\bibitem{dc_Co} F.~R.~Cohen,
\emph{The homology of $\mathcal{C}_{n+1}$-spaces, $n\ge 0$},
in: \emph{The homology of iterated loop spaces}, pp.~207--352, Lecture Notes in Math.,
vol.~533, Springer, Berlin, 1976; 
\MRh{0436146}

\bibitem{dc_CT} F. Cohen, L. Taylor, \emph{Computations of Gelfand-Fuks cohomology, the cohomology of function spaces, and the cohomology of configuration spaces}, in: \emph{Geometric applications of homotopy theory, I}, pp. 106--143, Lecture Notes in Math., vol. 657, Springer, Berlin, 1978; \MRh{0513543}.

\bibitem{Deligne} P. Deligne, \emph{Les immeubles des groupes de tresses g\'en\'eralis\'es}, Invent. Math. \textbf{17} (1972), 273--302; \MRh{0422673}.

\bibitem{dc_FN} E. Fadell, L. Neuwirth, \emph{Configuration spaces}, Math. Scand. \textbf{10} (1962), 111--118;  \MRh{0141126}.

\bibitem{dc_FR} M. Falk, R. Randell, 
\emph{The lower central series of a
fiber-type arrangement}, Invent. Math., \textbf{82} (1985), 77--88; 
\MRh{0808110}.

\bibitem{FaTC1} {M. Farber.} \emph{Topological complexity of motion
planning}, {Discrete Comput. Geom.} \textbf{29} (2003),
211--221; \MRh{1957228}.

\bibitem{dc_Fa03} M. Farber,  \emph{Instabilities of robot motion}, Topology Appl. \textbf{140} (2004), 245--266; 
\MRh{2074919}.

\bibitem{FarberGraphs1} M. Farber, \emph{Collision free motion planning on graphs}, in: \emph{Algorithmic Foundations of Robotics IV}, pp. 123--138, Springer Tracts in Adv. Robotics, vol. 17, Springer, Berlin, 2005.

\bibitem{dc_Fa05} M. Farber,  
{\em{Topology of robot motion planning}}, 
in: \emph{Morse Theoretic Methods in Non-linear 
Analysis and in Symplectic Topology}, pp. 185--230, 
NATO Science Series II: Mathematics, Physics and Chemistry, vol. 217, 
Springer, 2006; 
\MRh{2276952}.

\bibitem{FarberGraphs2} M. Farber,  \emph{Configuration spaces and robot motion planning algorithms}, in: \emph{Combinatorial and Toric Topology}, Lecture Note Series, Institute for Mathematical Sciences, National University of Singapore, World Scientific Publishing Co., 2017, to appear.

\bibitem{dc_FG} M. Farber, M. Grant, 
\emph{Topological complexity of configuration spaces}, Proc. Amer. Math. Soc. \textbf{137} (2009), 1841--1847; \MRh{2470845}.

\bibitem{dc_fgy} M. Farber, M. Grant, S. Yuzvinsky, 
\emph{Topological complexity of collision free motion planning algorithms in the presence of multiple moving obstacles}, in: \emph{Topology and Robotics}, pp. 75--83, Contemp. Math., vol. 438, Amer. Math. Soc., Providence, RI, 2007; \MRh{2359030}.

\bibitem{dc_FY} M.~Farber, S.~Yuzvinsky,
{\em{Topological robotics: Subspace arrangements and collision
free motion planning}}, in: {\em Geometry, topology, and 
mathematical physics}, pp.~145--156, 
Amer. Math. Soc. Transl. Ser. 2, vol.~212, 
Amer. Math. Soc., Providence, RI, 2004; 
\MRh{2070052}.

\bibitem{dc_NF} N. Fieldsteel, \emph{Topological complexity of graphic arrangements}, 
in: \emph{Topological Complexity and Related Topics}, Contemp. Math., 
Amer. Math. Soc., Providence, RI (this volume).

\bibitem{Ghrist} R. Ghrist, \emph{Configuration spaces and braid groups on graphs in robotics}, in: \emph{Knots, braids, and mapping class groups--papers dedicated to Joan S. Birman (New York, 1998)}, pp. 29--40, AMS/IP Stud. Adv. Math., vol. 24, Amer. Math. Soc., Providence, 2001; \MRh{1873106}.

\bibitem{GK} R. Ghrist, D. Koditschek, \emph{Safe cooperative robot dynamics on graphs}, SIAM J. Control Optim. \textbf{40} (2002), 1556--1575; \MRh{1882808}.

\bibitem{dc_GoGr} J. Gonz\'alez, M. Grant, \emph{Sequential motion planning of non-colliding particles in Euclidean spaces}, Proc. Amer. Math. Soc., \textbf{143} (2015), 4503--4512; \MRh{3373948}.

\bibitem{dc_GG} J. Gonz\'alez, B. Guti\'errez, \emph{Topological complexity of collision-free multi-tasking motion planning on orientable surfaces},  
in: \emph{Topological Complexity and Related Topics}, Contemp. Math., 
Amer. Math. Soc., Providence, RI (this volume).

\bibitem{ggy} J. Gonz\'alez, B. Guti\'errez, S. Yuzvinsky, \emph{The higher topological complexity of subcomplexes of products of spheres---and related polyhedral product spaces}, preprint 2015; \arxiv{1501.07474}.

\bibitem{Grant} M. Grant, \emph{Topological complexity, fibrations and symmetry}, Topology Appl. \textbf{159} (2012), 88--97; \MRh{2852952}.

\bibitem{dc_GLO} M. Grant, G. Lupton, J. Oprea, \emph{New lower bounds for the topological complexity of aspherical spaces}, Topology Appl. \textbf{189} (2015), 78--91; \MRh{3342573}.

\bibitem{GR-M} M. Grant, D. Recio-Mitter, \emph{Topological complexity of subgroups of Artin's braid groups},  
in: \emph{Topological Complexity and Related Topics}, Contemp. Math., 
Amer. Math. Soc., Providence, RI (this volume).

\bibitem{M2}  D. Grayson, M. Stillman, \emph{Macaulay2: a software system for research in algebraic geometry}, \\ available
at \href{http://www.math.uiuc.edu/Macaulay2}{\texttt{www.math.uiuc.edu/Macaulay2}}.

\bibitem{OS} P.~Orlik, L.~Solomon,
{\em Combinatorics and topology of complements of
hyperplanes}, Invent. Math. \textbf{56} (1980), 167--189; \MRh{0558866}.

\bibitem{OT}  P.~Orlik, H.~Terao,
{\em Arrangements of hyperplanes}, Grundlehren Math. Wiss.,
vol.~300, Springer-Verlag, New~York-Berlin-Heidelberg, 1992; \MRh{1217488}.

\bibitem{Ru} Y. Rudyak, \emph{On higher analogs of topological complexity}, Topology Appl. \textbf{157} (2010), 916--920;\\ \MRh{2593704}.

\bibitem{Scheirer} S. Scheirer, \emph{Topological complexity of $n$ points on a tree}, preprint 2017; \arxiv{1607.08185}.

\bibitem{dc_To} B. Totaro, \emph{Configuration spaces of algebraic varieties}, Topology \textbf{35} (1996), 1057--1067;\\ \MRh{1404924}.

\bibitem{dc_X} M. Xicot\'encatl, 
\emph{Orbit configuration spaces}, Thesis, University of Rochester, 1997.

\bibitem{YuzOS} S. Yuzvinsky, \emph{Orlik-Solomon algebras in algebra and topology}, Russian Math. Surveys \textbf{56} (2001), 293--364; \MRh{1859708}.

\bibitem{YuzHTC} S. Yuzvinsky, \emph{Higher topological complexity of Artin type groups}, 
in: {\em Configuration Spaces: Geometry, Topology and Representation Theory}, pp.~119--128, 
Springer INdAM Ser., vol. 14, Springer, Cham, 2016. 

\end{thebibliography}

\end{document}